\documentclass[a4paper]{article}
\pdfoutput=1

\usepackage{graphicx,wrapfig}
\usepackage{lipsum}
\usepackage{amsmath,amssymb,mathrsfs}
\usepackage{hyperref}
\usepackage{url}
\usepackage{mathtools}
\usepackage{changepage}
\usepackage{multirow}
\usepackage{wrapfig}
\usepackage{subfig}
\usepackage{color}
\usepackage{epsfig,enumerate,amsmath,amsfonts,amssymb,amsthm,mathrsfs,ifpdf}
\usepackage{indentfirst,relsize}
\usepackage{setspace,graphicx}
\usepackage{latexsym}
\usepackage[all]{xy}
\usepackage[usenames,dvipsnames]{pstricks}
\usepackage{pst-grad} 
\usepackage{pst-plot} 
\usepackage[margin = 2.50cm]{geometry}

\usepackage{algorithm}

\usepackage{algpseudocode}

\newcommand{\remove}[1]{}

\newtheorem{theo}{Theorem}
\newtheorem{lem}[theo]{Lemma}
\newtheorem{pre}[theo]{Proposition}
\newtheorem{coro}[theo]{Corollary}

\newcounter{Ca}[theo]
\newtheorem{ca}[Ca]{Case}
\newtheorem{defi}[theo]{Definition}
\newtheorem{remk}{Remark}

\title{Transitivity on subclasses of chordal graphs}

\author{Subhabrata Paul \and Kamal Santra }

\author{Subhabrata Paul\footnote{Department of Mathematics, IIT Patna, India} \and Kamal Santra\footnote{Department of Mathematics, IIT Patna, India} }

\date{}

\begin{document}

\maketitle
\begin{abstract}
	Let $G=(V, E)$ be a graph, where $V$ and $E$ are the vertex and edge sets, respectively. For two disjoint subsets $A$ and $B$ of $V$, we say $A$ \textit{dominates} $B$ if every vertex of $B$ is adjacent to at least one vertex of $A$ in $G$. A vertex partition $\pi = \{V_1, V_2, \ldots, V_k\}$ of $G$ is called a \emph{transitive $k$-partition} if $V_i$ dominates $V_j$ for all $i,j$, where $1\leq i<j\leq k$. The maximum integer $k$ for which the above partition exists is called \emph{transitivity} of $G$ and it is denoted by $Tr(G)$. The \textsc{Maximum Transitivity Problem} is to find a transitive partition of a given graph with the maximum number of partitions. It was known that the decision version of \textsc{Maximum Transitivity Problem} is NP-complete for chordal graphs [Iterated colorings of graphs, \emph{Discrete Mathematics}, 278, 2004]. In this paper, we first prove that this problem can be solved in linear time for \emph{split graphs} and for the \emph{complement of bipartite chain graphs}, two subclasses of chordal graphs. We also discuss Nordhaus-Gaddum type relations for transitivity and provide counterexamples for an open problem posed by J. T. Hedetniemi and S. T. Hedetniemi [The transitivity of a graph, \emph{J. Combin. Math. Combin. Comput}, 104, 2018]. Finally, we characterize transitively critical graphs having fixed transitivity.

\end{abstract}

{\bf Keywords.}
Transitivity, Split graphs, Complement of bipartite chain graphs, Nordhaus-Gaddum relations, Transitively critical graphs.

\section{Introduction}
Graph partitioning is one of the classical problems in graph theory. In a partitioning problem, the goal is to partition the vertex set (or edge set) into some parts with desired properties, such as independence, having minimum edges across partite sets, etc. In this article, we are interested in partitioning the vertex set into some parts such that the partite sets follow some domination relation among themselves. For a graph $G=(V,E)$, the \emph{neighbourhood} of a vertex $v\in V$ is the set of all adjacent vertices of $v$ and is denoted as $N_G(v)$. The \emph{degree} of a vertex $v$ in $G$, denoted as $\deg_G(v)$, is the number of edges incident to $v$. A vertex $v$ is said to \emph{dominate} itself and all its neighbouring vertices. A \emph{dominating set} of $G=(V,E)$ is a subset of vertices $D$ such that every vertex $x\in V\setminus D$ has a neighbour $y\in D$, that is, $x$ is dominated by some vertex $y$ of $D$. For two disjoint subsets $A$ and $B$ of $V$, we say $A$ \emph{dominates} $B$ if every vertex of $B$ is adjacent to at least one vertex of $A$. 

Graph partitioning problems, based on a domination relation among the partite sets, have been extensively studied in literature. Cockayne and Hedetniemi, in 1977, introduced the notion of \emph{domatic partition} of a graph $G=(V,E)$, where the vertex set is partitioned into $k$ parts, $\pi =\{V_1,V_2, \ldots, V_k\}$, such that each $V_i$ is a dominating set of $G$ \cite{cockayne1977towards}. The maximum order of such a domatic partition is called \emph{domatic number} of $G$ and it is denoted by $d(G)$. Another similar type of partitioning problem is \emph{the Grundy partition}. Christen and Selkow introduced a Grundy partition of a graph $G=(V,E)$ in 1979 \cite{CHRISTEN197949}. In the Grundy partitioning problem, the vertex set is partitioned into $k$ parts, $\pi =\{V_1,V_2, \ldots, V_k\}$, such that each $V_i$ is an independent set and for all $1\leq i< j\leq k$, $V_i$ dominates $V_j$. The maximum order of such a partition is called \emph{the Grundy number} of $G$ and it is denoted by $\Gamma(G)$. In 2004, Hedetniemi et al. introduced another such partitioning problem, namely \emph{upper iterated domination partition} \cite{erdos2003equality}. In an upper iterated domination partition, the vertex set is partitioned into $k$ parts, $\pi =\{V_1,V_2, \ldots, V_k\}$, such that for each $1\leq i\leq k$, $V_i$ is a minimal dominating set of $G\setminus (\cup_{j=1}^{i-1} V_j)$. The \emph{upper iterated domination number}, denoted by $\Gamma^*(G)$, is equal to the maximum order of such a vertex partition. Recently, in 2018, Haynes et al. generalized the idea of domatic partition and introduced the concept of \emph{upper domatic partition} of a graph $G$, where the vertex set is partitioned into $k$ parts, $\pi =\{V_1,V_2, \ldots, V_k\}$, such that for each $i, j$, with $1\leq i<j\leq k$, either $V_i$ dominates $V_j$ or $V_j$ dominates $V_i$ or both \cite{haynes2020upper}. The maximum order of such an upper domatic partition is called \emph{upper domatic number} of $G$ and it is denoted by $D(G)$. All these problems, domatic number \cite{chang1994domatic}, Grundy number \cite{hedetniemi1982linear, zaker2006results}, upper iterated number \cite{erdos2003equality}, upper domatic number \cite{haynes2020upper} have been extensively studied both from an algorithmic and structural point of view.

In this article, we study a similar graph partitioning problem, namely \emph{transitive partition}. In 2018, Hedetniemi et al. \cite{hedetniemi2018transitivity} have introduced this notion as a generalization of Grundy partition. A \emph{transitive $k$-partition} is defined as a partition of the vertex set into $k$ parts, $\pi =\{V_1,V_2, \ldots, V_k\}$, such that for all $1\leq i< j\leq k$, $V_i$ dominates $V_j$. The maximum order of such a transitive partition is called \emph{transitivity} of $G$ and is denoted by $Tr(G)$. The \textsc{Maximum Transitivity Problem (MTP)} is to find a transitive partition of a given graph with the maximum number of parts.
%
%
%
%
%
%
%
Note that a Grundy partition is a transitive partition with the additional restriction that each partite set must be independent. In a domatic partition $\pi =\{V_1,V_2, \ldots, V_k\}$ of $G$, since each partite set is a dominating set of $G$, we have domination property in both directions, that is, $V_i$ dominates $V_j$ and $V_j$ dominates $V_i$ for all $1\leq i< j\leq k$. However, in a transitive partition $\pi =\{V_1,V_2, \ldots, V_k\}$ of $G$, we have domination property in one direction, that is, $V_i$ dominates $V_j$ for $1\leq i< j\leq k$. In an upper domatic partition $\pi =\{V_1,V_2, \ldots, V_k\}$ of $G$, for all $1\leq i<j\leq k$, either $V_i$ dominates $V_j$ or $V_j$ dominates $V_i$ or both. The definition of each vertex partitioning problem ensures the following inequalities for any graph $G$. For any graph $G$, $1\leq \Gamma(G)\leq \Gamma^*(G)\leq Tr(G)\leq D(G)\leq n$.

In the introductory paper, J. T. Hedetniemi and S. T. Hedetniemi \cite{hedetniemi2018transitivity} showed, that the upper bound on the transitivity of a graph $G$ is $\Delta(G)+1$, where $\Delta(G)$ is the maximum degree of $G$. They also gave two characterizations for graphs with $Tr(G)=2$ and for graphs with $Tr(G)\geq 3$. They further showed that transitivity and Grundy number are the same for trees. Therefore, the linear-time algorithm for finding the Grundy number of a tree, presented in \cite{hedetniemi1982linear}, implies that we can find the transitivity of a tree in linear time as well. Also, for a subclass of bipartite graphs, namely bipartite chain graphs, MTP can be solved in linear time \cite{PaulSantra2022}. Moreover, for any graph, transitivity is equal to upper iterated domination number, that is, $\Gamma^*(G)=Tr(G)$ \cite{hedetniemi2018transitivity}, and the decision version of the upper iterated domination problem is known to be NP-complete for chordal graphs \cite{hedetniemi2004iterated}. Therefore, MTDP is NP-complete for chordal graphs as well. MTDP is also known to be NP-complete for perfect elimination bipartite graphs \cite{PaulSantra2022}. It is also known that every connected graph $G$ with $Tr(G)=k\geq 3$ has a transitive partition $\pi =\{V_1,V_2, \ldots, V_k\}$ such that $\lvert V_k \rvert$ = $\lvert V_{k-1} \rvert = 1$ and $\lvert V_{k-i} \rvert \leq 2^{i-1}$ for $2\leq i\leq k-2$ \cite{haynes2019transitivity}. This implies that MTP is fixed-parameter tractable \cite{haynes2019transitivity}. Also, graphs with transitivity at least $t$, for some integer $t$, have been characterized in \cite{PaulSantra2022}.

In this article, we study the computational complexity of the transitivity problem in subclasses of chordal graphs. The organization and main contributions of this article are summarized as follows. Section 2 contains basic definitions and notations that are followed throughout the article. Sections 3 and 4 describe two linear-time algorithms for split and for the complement of bipartite chain graphs, respectively. Section 5 deals with Nordhaus-Gaddum type relations for transitivity. In Section 6, we present a characterization of transitively vertex-edge critical graphs having fixed transitivity. Finally, Section 7 concludes the article.

%
%
%
%
%
%
%
%
%

\section{Notation and definition}

Let $G=(V, E)$ be a graph with $V$ and $E$ as its vertex and edge sets, respectively.  A graph $H=(V', E')$ is said to be a \emph{subgraph} of a graph $G=(V, E)$, if and only if $V'\subseteq V$ and $E'\subseteq E$. For a subset $S\subseteq V$, the \emph{induced subgraph} on $S$ of $G$ is defined as the  subgraph of $G$ whose vertex set is $S$ and edge set consists of all of the edges in $E$ that have both endpoints in $S$ and it is denoted by $G[S]$. The \emph{complement} of a graph $G=(V,E)$ is the graph $\overline{G}=(\overline{V}, \overline{E})$, such that $\overline{V}=V$ and $\overline{E}=\{uv| uv\notin E\}$. 


A subset of $S\subseteq V$, is said to be an \emph{independent set} of $G$, if every pair of vertices in $S$ are non-adjacent.  A subset of $K\subseteq V$, is said to be a \emph{clique} of $G$, if every pair of vertices in $K$ are adjacent. The cardinality of a clique of maximum size is called \emph{clique number} of $G$ and it is denoted by $\omega(G)$. A graph $G=(V, E)$ is said to be a \emph{split graph} if $V$ can be partitioned into an independent set $S$ and a clique $K$. 

A graph is called \emph{bipartite} if its vertex set can be partitioned into two independent sets. A bipartite graph $G=(X\cup Y,E)$ is called a \textit{bipartite chain graph} if there exists an ordering of vertices of $X$ and $Y$, say $\sigma_X= (x_1,x_2, \ldots ,x_{n_1})$ and $\sigma_Y=(y_1,y_2, \ldots ,y_{n_2})$, such that $N(x_{n_1})\subseteq N(x_{n_1-1})\subseteq \ldots \subseteq N(x_2)\subseteq N(x_1)$ and $N(y_{n_2})\subseteq N(y_{n_2-1})\subseteq \ldots \subseteq N(y_2)\subseteq N(y_1)$. Such ordering of $X$ and $Y$ is called a \emph{chain ordering} and it can be computed in linear time \cite{heggernes2007linear}. A graph $G$ is said to be a $2K_2$-free, if it does not contain a pair of independent edges as an induced subgraph. It is well-known that the class of bipartite chain graphs and $2K_2$-free bipartite graphs are the same. An edge between two non-consecutive vertices of a cycle is called a \emph{chord}. 
If every cycle in $G$ of length at least four has a chord, then $G$ is called a \emph{chordal graph}. 


\section{Transitivity in split graphs}
In this section, we design a linear-time algorithm for finding the transitivity of a given split graph. To design the algorithm, we first prove that the transitivity of a split graph $G$ can be either $\omega(G)$ or $\omega(G)+1$, where $\omega(G)$ is the size of a maximum clique in $G$. Further, we characterize the split graphs with the transitivity equal to $\omega(G)+1$.

\begin{lem}\label{SGTH1}
	Let $G=(S\cup K, E)$ be a split graph, where $S$ and $K$ are  an independent set and a clique of $G$, respectively. Also, assume that $K$ is the maximum clique of $G$, that is, $\omega(G)=|K|$. Then $\omega(G)\leq Tr(G)\leq \omega(G) +1$. Further, $Tr(G)=\omega(G)+1$ if and only if every vertex of $K$ has a neighbour in $S$. 
\end{lem}
\begin{proof}
	Note that $Tr(G)\geq \omega(G)$. As we can make a transitive partition $\pi=\{V_1, V_2, \ldots, V_{\omega(G)}\}$ of size $\omega(G)$ by considering each $V_i$ contains exactly one vertex from  maximum clique and all the other vertices in $V_1$. To prove that $Tr(G)\leq \omega(G) +1$, suppose $Tr(G)\geq \omega(G) +2$. Let $\pi=\{V_1, V_2, \ldots, V_ {\omega(G)+2}\}$ be a transitive partition of $G$. Since $\lvert K \rvert=\omega(G)$, there exist at least two sets in $\pi$, say $V_i$ and $V_j$ with $i<j$, such that $V_i$ and $V_j$ contains only vertices from $S$. Note that, in this case $V_i$ cannot dominate $V_j$ as $S$ is an independent set of $G$. Therefore, we have a contradiction. Hence, $\omega(G)\leq Tr(G)\leq \omega(G) +1$. 
	%
	%
	%
	%
	%
	
	Let every vertex of $K$ have a neighbour in $S$. Now consider a vertex partition of $G$, say $\pi=\{V_1, V_2, \ldots, V_{\omega(G)+1}\}$, such that $V_1 = S$ and for each $i>1$,  $V_i$ contains exactly one vertex from $K$. Since every vertex of $K$ has a neighbour in $S$, $V_1$ dominates every other partition in $\pi$. Moreover, as $K$ is a clique, each $V_i$, with $i>1$ dominates $V_j$ for all $2\leq i<j \leq \omega+1$. Hence, $\pi$ is a transitive partition of $G$. Now, since $Tr(G)\leq \omega(G) +1$, we have $Tr(G)=\omega(G)+1$.

	Conversely, let $Tr(G)=\omega(G)+1$ and $\pi=\{V_1,V_2, \ldots, V_{\omega(G) +1}\}$ be a transitive  partition of $G$. Note that if there exist two sets in $\pi$, that contain only vertices from $S$, then using similar arguments as before, we have a contradiction. Therefore, there exists exactly one set in $\pi$, say $V_l$, that contains only vertices from $S$ as $\omega(G)=|K|$. Hence, each set of $\pi$, except $V_l$, contains exactly one vertex from $K$. Suppose there exists a vertex, say $x$, in $K$ that has no neighbour in $S$ and also let $x\in V_p$ for some set $V_p$ in $\pi$. Note that there is no edge between the vertices of $V_p$ and $V_l$. Therefore, neither $V_p$ dominates $V_l$ nor $V_l$ dominates $V_p$. This contradicts the fact that $\pi$ is a transitive partition. Hence, every vertex of $K$ has a neighbour in $S$. 
\end{proof}

Based on the above lemma, we have the following algorithm for finding transitivity of a given split graph.

\begin{algorithm}[h]
	
	\caption{\textsc{Transitivity\_Split($G$)}}\label{Algo:Split}
	
	\begin{algorithmic}[1]
		
		\State  \textbf{Input:} A split graph $G=(V,E)$

		\State  \textbf{Output:} The transitivity of $G$, that is, $Tr(G)$
		
		\State Find a vertex partition of $V$ into $S$ and $K$, where $S$ and $K$ are an independent set and a clique of $G$, respectively and  $\omega(G)=|K|$.
		
		\ForAll {$v\in K$}
		
		\If {$v$ has no neighbour in $S$}
		
		\State $t=|K|$.
		
		\State {\bf break}
		
		\EndIf
		
		\EndFor
		
		\State $t=|K|+1$.
		
		\State \Return(t)
		
	\end{algorithmic}
	
\end{algorithm}

Note that the required vertex partition in line $3$ of Algorithm \ref{Algo:Split} can be computed in linear time \cite{hammer1981splittance}. Also, the for loop in line $4-7$ runs in O(n+m) time. Hence, we have the following theorem:

\begin{theo}
	The \textsc{Maximum Transitivity Problem} can be solved in linear time for split graphs.
\end{theo}

\section{Transitivity in the complement of bipartite chain graphs}

In this section, we find the transitivity of the complement of a bipartite chain graph, say $G$, by showing that the transitivity of $G$ is equal to the Grundy number of $G$. To do that, first, we show two essential properties, one of the transitive partitions and the second of the Grundy partitions for the complement of a bipartite chain graph.The proofs of these lemmas can be found in the Appendix.  

%

\begin{lem}\label{CCB1}
	Let $G=(X\cup Y, E)$ be the complement of a bipartite chain graph and $Tr(G)=k$. Then there exists a transitive partition $\pi = \{V_1, V_2, \ldots, V_k\}$ of $G$ such that either $\lvert V_1 \rvert=1$ or $V_1=\{x, y\}$, where $x\in X$ and $y\in Y$.	
\end{lem} 

\begin{proof}
	Let $\pi = \{V_1, V_2, \ldots, V_k\}$ be a transitive partition of $G$. If either $\lvert V_1 \rvert=1$ or $V_1=\{x, y\}$, we are done. Otherwise, exactly one of the following two cases happen:
	
	\begin{ca}\label{CA1}
		$|V_1|>2$ and $V_1$ contains vertices from $X$ and $Y$
	\end{ca}
	In this case, we can construct another transitive  partition $\pi'=\{V'_1, V'_2, \ldots, V'_k\}$ as follows: we retain only two vertices in $V'_1$, one from $X$, say $x$, and another from $Y$, say $y$. After that, we put every other vertex of $V_1$ into $V'_2$. Therefore, $V'_1=\{x, y\}$ for some $x\in X$ and $y\in Y$, $V'_2=V_2\cup(V_1\setminus\{x,y\})$ and $V'_i=V_i$, for all $i\geq 3$. Note that $V_1$ dominates every other set in $\pi'$ as both $X$ and $Y$ induces complete graph in $G$. Hence, $\pi'$ is a transitive partition of $G$.

	\begin{ca}\label{CA2}
		$|V_1| \geq 2$ and $V_1$ contains only vertices from $X$ or $Y$
	\end{ca}
	
	Without loss of generality, let us assume that $V_1=\{x_{t_1}, x_{t_2}, \ldots, x_{t_s}\}$ contains only vertices from $X=\{x_1, x_2, \ldots, x_{n_1}\}$. In this case, we show that there exists a vertex in $V_1$ which is adjacent with every vertex of $Y=\{y_1, y_2, \ldots, y_{n_2}\}$. Let $x_{t_i}$ be a vertex in $V_1$ which has maximum number of neighbours in $Y$. Let $Y_{t_i}=\{y_{p_{1}}, y_{p_2}, \ldots, y_{p_{s'}}\}$ be the neighbours of $x_{t_i}$ in $Y$. Suppose $|Y_{t_i}| < |Y|$. Consider any vertex from $(Y\setminus Y_{t_i})$, say $y$. Since $y$ is not adjacent to $x_{t_i}$ and $\pi$ is a transitive partition of $G$, $y$ must be adjacent with some vertex, say $x_{t_j}$, in $V_1$ other than $x_{t_i}$. Now, consider an arbitrary vertex $y' \in Y_{t_i}$. Note that $x_{t_j} y' \in E$ because otherwise $\{x_{t_i}, y, x_{t_j}, y'\}$ induces a $2K_2$ in $\overline{G}$ (see Fig. \ref{fig:complement_chain}). This contradicts the fact that $\overline{G}$ is a bipartite chain graph. Now, since $y'$ is chosen arbitrarily from $Y_{t_i}$, we can say that $x_{t_j}$ is adjacent to every vertex of $Y_{t_i}$. This contradicts the fact that $x_{t_i}$ is a vertex in $V_1$ which has maximum number of neighbours in $Y$. Hence, there exists a vertex, say $x$, in $V_1$ which is adjacent with every vertex of $Y$. We consider the partition $\pi'=\{V'_1, V'_2, \ldots, V'_k\}$, where $V'_1=\{x\}$, $V'_2=V_2\cup (V_1\setminus \{x\})$ and $V'_i=V_i$, for all $i\geq 3$. Clearly, $\pi'$ is a transitive partition of $G$.

	\begin{figure}[htbp!]
		\centering
		\includegraphics[scale=0.42]{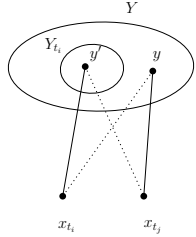}
		\caption{$2K_2$ in $\overline{G}$}
		\label{fig:complement_chain}
		
	\end{figure}

Hence, we can always construct a transitive partition $\pi = \{V_1, V_2, \ldots, V_k\}$ of $G$ such that either $\lvert V_1 \rvert=1$ or $V_1=\{x, y\}$, where $x\in X$ and $y\in Y$. 

\end{proof}

\begin{lem}\label{CCB2}
	Let $G=(X\cup Y, E)$ be the complement of a bipartite chain graph. Also, let $\pi = \{V_1, V_2, \ldots, V_k\}$ be a Grundy partition of $G$ with $\Gamma(G)=k$ and $X'_G=\{x\in X| x\in V_i \; \text{and} \;\lvert V_i \rvert=1\}$ and $Y'_G=\{y\in Y| y\in V_i \; \text{and} \;\lvert V_i \rvert=1\}$. Then exactly one of the following two cases is true:
	
	(i) Both $\lvert X'_G \rvert$ and $\lvert Y'_G \rvert$ cannot be empty simultaneously.
	
	(ii) The graph $G$ is the disjoint union of $K_{|X|}$ and $K_{|Y|}$ and $|X|=|Y|$.
	
\end{lem}

\begin{proof}
	
		As $\pi$ is a Grundy partition of  the complement of a bipartite graph, then every set in $\pi$ is either size two or size one. If $|X|\neq |Y|$, then we always have either $\lvert X'_G \rvert \not=\phi$ or  $\lvert Y'_G \rvert \not= \phi$. So, we assume that $|X|=|Y|$. Now, if either $\lvert X'_G \rvert \not=\phi$ or  $\lvert Y'_G \rvert \not= \phi$, then we are done. Hence, let us assume that $\lvert X'_G \rvert =\phi$ and  $\lvert Y'_G \rvert = \phi$. Note that, in this case, each set in $\pi$ contains exactly one vertex from $X$ and  another vertex from $Y$, because each set has to be an independent set. Let for all $1\leq i\leq k$, $V_i=\{x_i, y_i\}$. Consider the following cases:

	\begin{ca}\label{one_edge_present_1}
		There exists $V_p, V_q\in \pi$ such that $x_py_q \in E$ and $x_qy_p\in E$
	\end{ca}
	In this case, $\{x_p,y_p,x_q,y_q\}$ forms a $2K_2$ in $\overline{G}$, since $V_p$ and $V_q$ are independent sets. This contradicts the fact that $\overline{G}$ is a bipartite chain graph.
	
	\begin{ca}\label{one_edge_present_2}
		There exists $V_p, V_q\in \pi$ such that $x_py_q \in E$ and $x_qy_p\notin E$
	\end{ca}
	Without loss of generality, let $p<q$. In this case, we form a new partition $\pi'=\{V'_1, V'_2, \ldots, V'_{k+1}\}$ as follows: $V'_i=V_i$ for all $1\leq i\leq (q-1)$ and $i\neq p$, $V'_p=\{x_q, y_p\}$, $V'_j=V_{j+1}$ for $q\leq j\leq (k-1)$, $V'_k =\{x_p\}$ and $V'_{k+1}=\{y_q\}$. Note that for all $1\leq i\leq (k-1)$, as $V'_i$ contains both vertices from $X$ and $Y$, clearly, $V'_i$ dominates every other partition in $\pi'$. Also, since $x_py_q \in E$, $V'_k$ dominates $V'_{k+1}$. Therefore, $\pi'$ is a Grundy partition which is a contradiction to the fact that $\Gamma(G)=k$.
	
	\begin{ca}\label{other_edge_present_3}
		There exists $V_p, V_q\in \pi$ such that $x_py_q \notin E$ and $x_qy_p\in E$
	\end{ca}
	Similar to the previous case.
	
	\begin{ca}\label{other_edge_present_4}
		For every pair of sets $V_p, V_q\in \pi$, $x_py_q \notin E$ and $x_qy_p\notin E$
	\end{ca}
	In this case, there is no edge between $X$ and $Y$ in $G$. Hence, $G$ is the disjoint union of $K_{|X|}$ and $K_{|Y|}$.

	Hence, we have the above lemma. 
	
\end{proof}

Now, we are ready to show that the transitivity and the Grundy number are equal for the complement of a bipartite chain graph.

\begin{theo}\label{CCB3}
	Let $G=(X\cup Y, E)$ be the complement of a bipartite chain graph. Then $\Gamma(G)=Tr(G)$.
\end{theo}

\begin{proof}
	We use induction on $n$, where $n$ is the number of vertices of $G$. If $n=1$, then $\Gamma(G)=Tr(G)=1$ trivially. For $n=2$, $G$ is either $K_2$ or $\overline{K}_2$ and therefore, $\Gamma(G)=Tr(G)$. Let us assume that the induction hypothesis is true, that is, $\Gamma(G)=Tr(G)$ for the complement of every bipartite chain graph having less than $n$ vertices. Let us consider a transitive partition $\pi = \{V_1, V_2, \ldots, V_k\}$ of $G$ with $Tr(G)=k$. By Lemma \ref{CCB1}, we can assume that $\lvert V_1 \rvert =1$ or $V_1=\{x, y\}$ for some $x\in X$ and $y\in Y$. Let $H= G \setminus V_1$. Note that $H$ is also the complement of a bipartite chain graph, since deleting a vertex from $X$ (or $Y$) does not change the chain ordering of the remaining vertices. By induction hypothesis, we have $\Gamma(H)=Tr(H)$. Moreover, note that $Tr(H)=k-1$. Hence, we have $\Gamma(H)=Tr(H)= k-1$. Let $\pi'= \{V'_1, V'_2, \ldots, V'_{k-1}\}$ be a Grundy partition of $H$. Now, if $V_1=\{x\}$ (or \{y\}), then $x$ (correspondingly $y$) is adjacent to every vertex of $G$ because $\pi$ is a transitive partition of $G$. Therefore, $\pi''= \{V_1, V'_1, V'_2, \ldots, V'_{k-1}\}$ forms a Grundy partition of $G$ which implies $\Gamma(G)\geq k =Tr(G)$. Also, for any graph we know that $\Gamma(G)\leq Tr(G)$, hence $\Gamma(G)=Tr(G)$. So, let us assume that $V_1=\{x, y\}$ for some $x\in X$ and $y\in Y$. Now, if $xy\notin E$, that is, $V_1$ is an independent set, then from a Grundy partition of $H$ of order $(k-1)$ we can construct a Grundy partition of $G$ of order $k$ by appending $V_1$. Then by similar argument, we have $\Gamma(G)=Tr(G)$. So, we assume that $xy\in E$. Since, $H=G\setminus V_1$ is the complement of a bipartite chain graph, by induction hypothesis $\Gamma(H)=Tr(H)$. Now, by Lemma \ref{CCB2}, we can assume that $H$ has a Grundy partition, say $\pi' = \{V'_1, V'_2, \ldots, V'_{k-1}\}$, such that either $|X'_H| \not=\phi$ or $|Y'_H|\not=\phi$, where $X'_H$ and $Y'_H$ is defined in a similar way as in Lemma \ref{CCB2}  or $H$ is the disjoint union of $K_{|X_H|}$ and $K_{|Y_H|}$ and $|X_H|=|Y_H|$, where $X_H=X\cap H$ and $Y_H=Y\cap H$. If $H$ is the disjoint union of $K_{|X_H|}$ and $K_{|Y_H|}$ and $|X_H|=|Y_H|$, then $\pi''=\{V'_1, V'_2, \ldots, V'_{k-1}, \{x\}, \{y\}\}$ forms a transitive partition of $G$ of order $(k+1)$. This is a contradiction to the fact that $Tr(G)=k$. So, we assume that $H$ has a Grundy partition, say $\pi' = \{V'_1, V'_2, \ldots, V'_{k-1}\}$, such that either $|X'_H| \not=\phi$ or $|Y'_H|\not=\phi$. The remaining proof is done by dividing into the following four cases:


	\begin{ca} \label{CAA1}
		Every vertex of $X'_H$ and $Y'_H$ are adjacent to $y$ and $x$, respectively
	\end{ca}
	
	In this case, consider the vertex partition $\pi''=\{V'_1, V'_2, \ldots, V'_{k-1}, \{x\}, \{y\}\}$. Clearly, $\pi''$ is a transitive partition of $G$ of order $(k+1)$. This is a contradiction to the fact that $Tr(G)=k$.
	
	\begin{ca} \label{CAA2}
		Every vertex of $X'_H$ is adjacent to $y$ but there exists a vertex $y_t\in Y'_H$ such that $xy_t\notin E$
	\end{ca}
	Let $y_t\in V'_p$. In this case, let us consider the vertex partition $\pi''=\{U_1, U_2, \ldots, U_k\}$, where $U_1=\{x, y_t\}$, $U_i=V'_{i-1}$, for $2\leq i\leq p$, $U_{p+1}=\{y\}$ and $U_j=V'_{j-1}$, for all $p+2\leq j \leq k$. Clearly, the partition $\pi''$ forms a Grundy partition of $G$ which implies $\Gamma(G)\geq k =Tr(G)$. Since, for any graph $\Gamma(G)\leq Tr(G)$, therefore, $\Gamma(G)=Tr(G)$.
	
	\begin{ca} \label{CAA3}
		There exists a vertex $x_s\in X'_H$ such that $yx_s\notin E$ but every vertex of $Y'_H$ is adjacent to $x$
	\end{ca}
	This case is similar to Case \ref{CAA2}.
	
	\begin{ca} \label{CAA4}
		There exists a vertex $x_s\in X'_H$ such that $yx_s\notin E$ and there exists a vertex $y_t\in Y'_H$ such that $xy_t\notin E$
	\end{ca}
	In this case, $\{x, y_t, y, x_s\}$ induces a $2K_2$ in $\overline{G}$. This is a contradiction to the fact that $\overline{G}$ is a bipartite chain graph.
	
	Hence, for the complement of a bipartite chain graph $G$, $\Gamma(G)=Tr(G)$.
\end{proof}

It was proved in \cite{zaker2006results} that for the complement of a bipartite graph, $\Gamma(G)=n-p$, where $n$ is the number of vertices of $G$ and $p$ is the cardinality of a minimum edge dominating set of $\overline{G}$. We also know that the minimum edge dominating set of a bipartite chain graph can be computed in linear time \cite{verma2020grundy}. Therefore, we have the following corollary:

\begin{coro}
	The transitivity of the complement of bipartite chain graphs can be computed in linear time.
\end{coro}

\begin{remk}
	Identifying graphs with equal transitivity and Grundy number was posed as an open question in \cite{hedetniemi2018transitivity}. Theorem \ref{CCB3} partially answers this question by showing that the complement of bipartite chain graphs form such a graph class.
\end{remk}

\section{Nordhaus–Gaddum type bounds for transitivity}
Let $G=(V,E)$ be  a simple graph. A proper $k$-coloring of $G$ is a function $c$ from $V$ to $\{1, 2, \ldots, k\}$ such that $c(u)\neq c(v)$ if and only if $uv\in E$. The minimum value of $k$ for which a proper coloring exists is called \emph{chromatic number} of $G$ and it is denoted by $\chi(G)$. In 1956, Nordhaus and Gaddum \cite{nordhaus1956complementary} studied the chromatic number of a graph $G$ and its complement $\overline{G}$. They established lower and upper bound for the product and the sum of $\chi(G)$ and $\chi(\overline{G})$ in terms of the number of vertices of $G$. Since then, any bound on the sum or the product of a parameter of a graph $G$ and its complement $\overline{G}$ is called a Nordhaus–Gaddum type inequality. In this section, we study Nordhaus and Gaddum type relations for transitivity.


From \cite{PaulSantra2022}, it is known that for a bipartite chain graph $G$, $Tr(G)=t+1$, where $t$ is the maximum integer such that $G$ contains either $K_{t,t}$ or  $K_{t,t}-\{e\}$ as an induced subgraph. Let $\sigma_X= (x_1,x_2, \ldots ,x_{n_1})$ and $\sigma_Y=(y_1,y_2, \ldots ,y_{n_2})$ be the chain ordering of $G$. Because of this chain ordering if $x_py_p\in E$ for some $p$, then $\{x_1,x_2, \ldots ,x_{p}\}$ and $\{y_1,y_2, \ldots ,y_{p}\}$ induces a complete bipartite graph. Therefore, it follows that if $j$ is the maximum index such that $x_jy_j\in E$, then
$$Tr(G) = \begin{cases}
	j+2  & x_{j+1}y_{j}, x_jy_{j+1}\in E \\
	j+1 &  \text{otherwise}
\end{cases}$$
We know that for the complement of a bipartite graph $G$, $\Gamma(G)=n-p$, where $n$ is the number of vertices of $G$ and $p$ is the cardinality of a minimum edge dominating set of $\overline{G}$ \cite{zaker2006results}. Therefore, from Theorem \ref{CCB3} we have, for a bipartite chain graph $G$, $Tr(\overline{G})=n-p$, where $n$ is the number of vertices of $G$ and $p$ is the cardinality of a minimum edge dominating set of $G$. Also from \cite{verma2020grundy}, we know that for a bipartite chain graph, $p$ is equal to the maximum index $j$ such that $x_jy_j\in E$. Therefore, we have the following Nordhaus-Gaddum relation for transitivity in bipartite chain graphs:

\begin{theo}\label{theo:NGbipartitechain}
	Let $G=(X\cup Y, E)$ be a bipartite chain graph with the ordering $\sigma_X= (x_1,x_2, \ldots ,x_{n_1})$ and $\sigma_Y=(y_1,y_2, \ldots ,y_{n_2})$ as its chain ordering, that is, $N(x_{n_1})\subseteq N(x_{n_1-1})\subseteq \ldots \subseteq N(x_2)\subseteq N(x_1)$ and $N(y_{n_2})\subseteq N(y_{n_2-1})\subseteq \ldots \subseteq N(y_2)\subseteq N(y_1)$. Also assume that $j$ is the maximum index such that $x_jy_j\in E$. Then, 
	$$Tr(G) + Tr(\overline{G})= \begin{cases}
		n+2  & x_{j+1}y_{j}, x_jy_{j+1}\in E \\
		n+1 &  \text{otherwise}
	\end{cases}$$
\end{theo}

Based on Lemma \ref{SGTH1} and the fact that the complement of a split graph is also a split graph, we have the following Nordhaus-Gaddum relation for transitivity in split graphs, whose proof can be found in Appendix. 

\begin{theo}\label{theo:NGSplit}
	Let $G=(S\cup K, E)$ be a split graph, where $S$ and $K$ are the independent set and clique of $G$, respectively. Also, assume that $K$ is the maximum clique of $G$, that is, $\omega(G)=|K|$. Then, 
	$$Tr(G) + Tr(\overline{G})= \begin{cases}
		n+2  & \text{if, in $G$, every vertex of K has a neighbour in S} \\
		n+1 &  \text{otherwise}
	\end{cases}$$
\end{theo}
\begin{proof}
	
	Let us assume that, in $G$, every vertex of $K$ has a neighbour in $S$. Then by Lemma \ref{SGTH1}, $Tr(G)= |K|+1$. Note that, in $\overline{G}$, $S$ forms a clique. Moreover, none of the vertices in $K$ is adjacent to every vertex of $S$ as, in $G$,  every vertex of $K$ has a neighbour in $S$. Therefore, $\omega(\overline{G})= |S|$. Also, note that, in $\overline{G}$, every vertex of $S$ has a neighbour in $K$. If not, let us assume that there is a vertex in $\overline{G}$, say $x\in S$, which has no neighbour in $K$. Then $(K\cup \{x\})$ forms a clique in $G$, which contradicts the fact that $\omega(G)=|K|$. Therefore, the complement of $G$, that is, $\overline{G}=(K\cup S, \overline{E})$ is a split graph, where $K$ is an independent set and $S$ is a clique with $\omega(\overline{G})=|S|$. Also, in $\overline{G}$, every vertex of $S$ has a neighbour in $K$. Hence by Lemma \ref{SGTH1}, $Tr(\overline{G})= |S|+1$. This implies that, in this case, 
	$$Tr(G) + Tr(\overline{G})= |K|+1+|S|+1= n+2.$$
	
	On the other hand, let us assume that, in $G$, there exists a vertex, say $x$, in $K$ which has no neighbour in $S$. Then by Lemma \ref{SGTH1}, $Tr(G)= |K|$. Let $S'=(S\cup \{x\})$ and $K'= K\setminus \{x\}$. Note that, in $\overline{G}$, $S'$ forms a maximum clique in $\overline{G}$, that is, $\omega(\overline{G})=|S'|$. Also, since $K$ forms an independent set in $\overline{G}$, $x$ is not adjacent to any vertex of $K'$. Therefore, the complement of $G$, that is, $\overline{G}=(K'\cup S', \overline{E})$ is a split graph, where $K'$ is an independent set and $S'$ is a clique with $\omega(\overline{G})=|S'|$. Also, in $\overline{G}$, there exists a vertex, namely $x$, in $S'$ which has no neighbour in $K'$. Hence by Lemma \ref{SGTH1}, $Tr(\overline{G})= |S'|$. This implies that, in this case, 
	$$Tr(G) + Tr(\overline{G})= |K|+|S'|= |K|+|S|+1=n+1.$$
	
	Hence, the Nordhaus-Gaddum relation holds for split graphs. 
	
\end{proof}

\begin{remk}
	In \cite{hedetniemi2018transitivity}, Hedetniemi and Hedetniemi posed the following open question about the sum of $Tr(G)$ and $Tr(\overline{G})$: for any graph $G$, is $Tr(G)+Tr(\overline{G})=n+1$ if and only if $G=K_n$ or $G=\overline{K}_n$? Theorem \ref{theo:NGbipartitechain} and \ref{theo:NGSplit} show the existence of some bipartite chain graph and split graph, respectively, for which $Tr(G)+Tr(\overline{G})=n+1$. Moreover, also for $K_{n,n}$, $Tr(G)+Tr(\overline{G})=2n+1$ which shows another counter example for the above mentioned open question.
\end{remk}

\section{Transitively critical graphs}

The concept of transitively critical graph was introduced by Haynes et al. in \cite{haynes2019transitivity}. A graph $G=(V, E)$ is said to be  \emph{transitively vertex critical} (\emph{transitively edge critical}) if deleting any vertex from $V$ (respectively, edge from $E$) results in a graph whose transitivity is less than $Tr(G)$. A transitively vertex critical (transitively edge critical) graph $G$ with $Tr(G)=k$ is called by $Tr_k^v$-critical (respectively, $Tr_k^{e}$-critical). Characterizations of vertex critical graph have been studied in \cite{haynes2019transitivity}for some small values of $k$. In this section, we introduce a generalization of transitively critical graphs, namely \emph{transitively vertex-edge critical} graphs and give characterization of such graphs for every fixed value of $k$. Using this characterization, we then characterize transitively edge critical graphs for every fixed value of $k$.

A transitively vertex-edge critical graph is basically a graph which is both transitively vertex and edge critical. The formal definition is as follows:

\begin{defi}
	A graph $G=(V, E)$ is said to be a transitively vertex edge-critical graph if deleting any element from $V\cup E$ results in a graph whose transitivity is less than $Tr(G)$. A transitively vertex edge-critical graph $G$ with $Tr(G)=k$ is called $Tr_k^{(v,e)}$-critical.
\end{defi}

Note that unlike transitively edge critical graphs, every transitively vertex-edge critical graph is connected. The graph ${K}_1$ is the only connected graph with $Tr(G)=1$ and it is both transitively edge and vertex critical. Therefore, the only $Tr_1^{(v,e)}$-critical is $K_1$.

The following proposition characterizes the $Tr_2^{(v,e)}$-critical graphs.

\begin{pre}
	The only $Tr_2^{(v,e)}$-critical graph is $K_2$.
\end{pre}
\begin{proof}\label{ve-critical_2}
	Clearly, transitivity of $K_2$ is $2$. Also, if we remove any edge or any vertex from $K_2$, we are left with $\overline{K}_2$ or $K_1$ and  transitivity of those graphs is $1$. Hence, $K_2$ is a  $Tr_2^{(v,e)}$-critical.
	
	Let $G$ be a $Tr_2^{(v,e)}$-critical graph with $n$ vertices. Since $Tr(G)=2$, $G$ is a disjoint union of stars which is shown by Hedetniemi et al. \cite{hedetniemi2018transitivity}. As $G$ is a vertex critical graph, which implies $G$ must be a connected graph \cite{haynes2019transitivity}. Therefore, $G$ is a star. Now, if $G$ contains more than one edge, then removal of that edge from $G$ does not decrease the transitivity, which contradicts the fact that $G$ is a transitively edge-critical. Therefore, $G$ can only be $K_2$. Hence, the only $Tr_2^{(v,e)}$-critical graph is $K_2$. 
\end{proof}

Next, we generalize the characterization for $Tr_k^{(v, e)}$-critical graphs for $k\geq 3$. To this end, we recall the concept of $t$-atom, which was introduced by Zaker in \cite{zaker2006results}. For the sake of completeness, we give the definition of $t$-atom here.

\begin{defi}[\cite{zaker2006results}]
	A $t$-atom is defined in a recursive way as follows:
	\begin{enumerate}
		
		\item The only $1$-atom is $K_1$.
		\item Let $H=(V,E)$ be any $(t-1)$-atom with $n$ vertices. Consider an independent set ${I_r}$ on $r$ vertices for any $r\in \{1,2,\ldots n\}$. For that fixed $r$, consider a $r$ vertex subset $W$ of $V$ and add a perfect matching between the vertices of ${I_r}$  and $W$. Then join an edge between each vertex of $V\setminus W$ and an (and to only one) arbitrary vertex of ${I_r}$. The resultant graph $G$ is a $t$-atom.
	\end{enumerate}
	
\end{defi}

The set of $t$-atoms is denoted by $\mathcal{A}_t$. The following lemma describes the transitively vertex edge-critical graph with transitivity $k$.

\begin{lem}\label{k_atom-vertex-edge-critical}
	If $G$ is  $Tr_k^{(v,e)}$-critical, then $G\in \mathcal{A}_k$.
	
\end{lem}

\begin{proof}
	
	For an integer $t$, $Tr(G)\geq t$ if and only if $G$ contains a $t$-atom as a subgraph, which is shown by Paul and Santra \cite{PaulSantra2022}.
	Since the transitivity of $G$ is $k$, therefore, $G$ contains a $k$-atom as a subgraph. Let $H\in \mathcal{A}_k$ and $G$ contains $H$ as a subgraph. Since the $Tr(H)\geq k$ and $Tr(G)\geq Tr(H)$, then $Tr(H)=k$. If $G$ has an edge other than edges of $H$, then removal of that edge from $G$ does not decrease the transitivity, which contradicts the fact that $G$ is a transitively edge-critical. Also, $G$ cannot contain more vertex than $H$, as $G$ is a transitively vertex critical graph too. Therefore, $G= H$. 
\end{proof}

The only $3$-atoms are $K_3$ and $P_4$. Also the graphs $K_3$ and $P_4$ are both $Tr_3^{v}$-critical and $Tr_3^{e}$-critical. Therefore, the converse of Lemma \ref{k_atom-vertex-edge-critical} is true for $k=3$. Hence, we have the following corollary.
\begin{coro}
	The only $Tr_3^{(v,e)}$-critical graphs are $K_3$ or $P_4$.
\end{coro}

\begin{figure}[htbp!]
	\centering
	\includegraphics[scale=0.65]{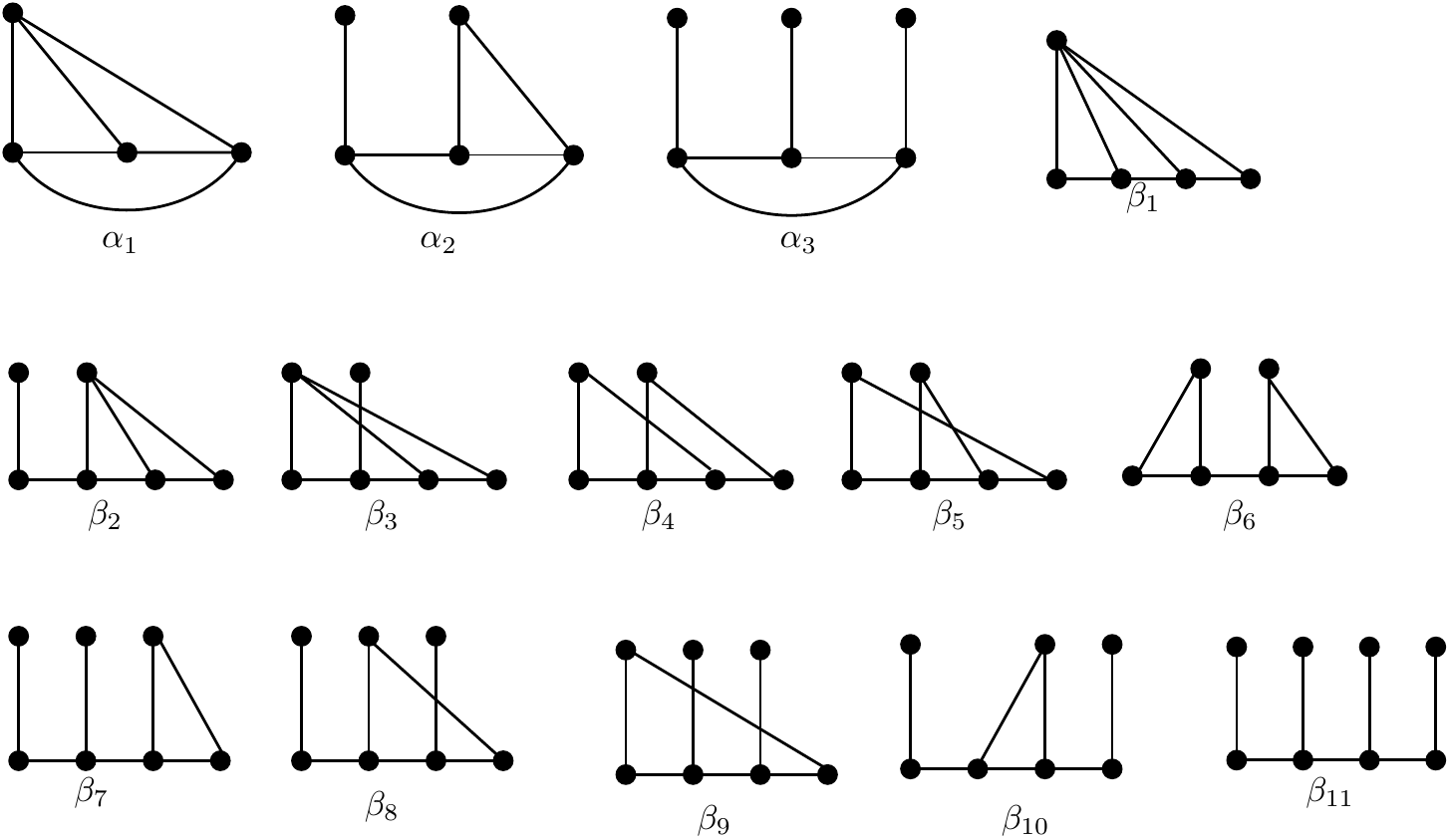}
	\caption{The class $\mathcal{A}_4$.}
	\label{fig:lits_4_atom}
\end{figure}

For $k=4$, the converse of Lemma \ref{k_atom-vertex-edge-critical} is not true. The class of graphs $\mathcal{A}_4$ is illustrated in Figure \ref{fig:lits_4_atom}. Note that every graph in $\mathcal{A}_4$, has transitivity equal to $4$ but only $\beta_2$ is not transitively edge-critical. Therefore, we have the following corollary.

\begin{coro}\label{Transitively 4-vertex-edge-critical1 }
	A graph	$G$ is $Tr_4^{(v,e)}$-critical if and only if $G\in \mathcal{A}_4'= (\mathcal{A}_4 \setminus \{\beta_2\})$.
\end{coro}

Generalizing this result, we have the following main theorem.

\begin{theo}
	Let $\mathcal{A}_k$ be the set of all $k$-atoms and $\mathcal{B}_k$ be the set of $k$-atoms which are neither $Tr_k^{(v,e)}$-critical nor have transitivity equal to $k$. A graph $G$ is $Tr_k^{(v,e)}$-critical if and only if $G\in \mathcal{A}_k'= (\mathcal{A}_k \setminus \mathcal{B}_k)$.
\end{theo}

\begin{proof}
	Let $G$ be a $Tr_k^{(v,e)}$-critical. Since the transitivity of $G$ is $k$, then $G$ contains a $k$-atom as a subgraph, as for an integer $t$, $Tr(G)\geq t$ if and only if $G$ contains a $t$-atom as a subgraph, which is shown by Paul and Santra \cite{PaulSantra2022}. Let $H$ be a $k$-atom and $G$ contains $H$ as a subgraph. Since the $Tr(H)\geq k$ and $Tr(G)\geq Tr(H)$, then $Tr(H)=k$. If $G$ has an edge other than edges of $H$, then removal of that edge from $G$ does not decrease the transitivity, which contradicts the fact that $G$ is a transitively edge critical. Also, $G$ cannot contain more vertex than $H$, as $G$ is also a transitively vertex critical graph. Hence, $G\in \mathcal{A}_k'= \mathcal{A}_k \setminus \mathcal{B}_k$.

\end{proof}

Next, we characterize the $Tr_k^{e}$-critical graphs for a fixed value of $k$. For this characterization, we first show the following relation between $Tr_k^{(v,e)}$-critical and $Tr_k^{e}$-critical graphs. The proof can be found in the Appendix.

\begin{theo}\label{relation_between_edge-crtical and ve-critical}
	A graph $G$ with $n$ vertices and $Tr(G)=k$ is $Tr_k^e$-critical if and only if $G=H\cup \overline{K}_{n-n_H}$, where $H$ is a $Tr_k^{(v,e)}$-critical graph having $n_H$ vertices.
\end{theo}

\begin{proof}
	
	Let $G=H\cup \overline{K}_{n-n_H}$, where $H$ is a $Tr_k^{(v,e)}$-critical graph having $n_H$ vertices. Since $H$ is a $Tr_k^{(v,e)}$-critical, then removal of any edge from H decreases its transitivity. Therefore, $G$ is a transitively edge-critical graph.
	
	Conversely, let $G$ be a $Tr_k^e$-critical graph. Also, let $G=C_1\cup C_2\cup \ldots \cup C_t$, where $C_i$ are the connected components of $G$. Since $Tr(G)=\max\{Tr(C_i) | 1\leq i \leq t\}$, then we may assume $Tr(G)=Tr(C_1)$. If there is  an edge in $C_i$, for any $i\geq 2$, then removal of that edge from $G$ does not decrease the transitivity, which contradicts the fact that $G$ is a transitively edge critical. Therefore, $C_i=K_1$, for all $i\geq 2$. Next, we prove that $C_1$ is a $Tr_k^{(v,e)}$-critical. Clearly, $C_1$ is a transitively edge critical as $G$ is a transitively edge critical. If $C_1$ is not a transitively vertex critical, then there exists a vertex in $C_1$, say $x$ such that $Tr(C_1)=Tr(C_1\setminus \{x\})=k$. Let $e$ be an edge incident to $x$ in $C_1$ and consider the graph $G'=C_1\setminus \{x\}$ and $G''=C_1\setminus \{e\}$. As, $G'\subseteq G''\subseteq C_1$, so $k=Tr(G')\leq Tr(G'')\leq Tr(C_1)=k$. Therefore, $Tr(G'')=k$, which contradicts the fact that $G$ is a transitively edge critical. So, $C_1$ is a transitively vertex critical graph. Hence, $C_1$ is a $Tr_k^{(v,e)}$-critical graph. Therefore, $G=C_1\cup \overline{K}_{n-n_{C_1}}$. 
\end{proof}

The characterization of $Tr_k^{e}$-critical graphs follows immediately from the above theorem.

%
%

\begin{coro}
	Let $\mathcal{A}_k$ be the set of $k$-atoms and $\mathcal{B}_k$ be the set of $k$-atoms which are neither $Tr_k^{(v,e)}$-critical nor have transitivity equal to $k$. A graph $G$ with $n$ vertices, is $Tr_k^{e}$-critical if and only if $G=H\cup \overline{K}_{n-n_H}$, where $H\in \mathcal{A}_k'= \mathcal{A}_k \setminus \mathcal{B}_k$.
\end{coro}

\section{Conclusion}
In this paper, we have proved that the transitivity of a given split and the complement of bipartite chain graphs can be computed in linear time. Then, we have discussed Nordhaus-Gaddum type relations for transitivity in split graphs and bipartite chain graphs and have given counter-examples to an open question posed in \cite{hedetniemi2018transitivity}. We have also studied transitively vertex-edge critical graphs. It would be interesting to investigate the complexity status of this problem in other subclasses of chordal graphs. Designing an approximation algorithm for this problem would be another challenging open problem.

\section*{Acknowledgements:} Subhabrata Paul is supported by the SERB MATRICS Research Grant (No. MTR/2019/000528). The work of Kamal Santra is supported by the Department of Science and Technology (DST) (INSPIRE Fellowship, Ref No: DST/INSPIRE/ 03/2016/000291), Govt. of India.

\bibliographystyle{alpha}
\bibliography{TC_Bib}

\end{document}